\documentclass{amsart}

\newcommand{\Z}{\ensuremath{\mathbf Z}}

\newcommand{\R}{\ensuremath{\mathbf R}}
\newtheorem{theorem}{Theorem}
\newcommand{\bt}{\begin{theorem}}
\newcommand{\et}{\end{theorem}}
\newtheorem{lemma}{Lemma}
\newcommand{\bl}{\begin{lemma}}
\newcommand{\el}{\end{lemma}}
\newtheorem{example}{Example}
\newcommand{\bex}{\begin{example}}
\newcommand{\eex}{\end{example}}
\newcommand{\beq}{\begin{equation}}
\newcommand{\eeq}{\end{equation}}
\newcommand{\benum}{\begin{enumerate}}
\newcommand{\eenum}{\end{enumerate}}

\title[Inverse problems for linear forms]{Inverse problems for linear forms over finite sets of integers}

\author{Melvyn B. Nathanson}

\address{Department of Mathematics\\Lehman College (CUNY)\\Bronx, New York 10468}
\email{melvyn.nathanson@lehman.cuny.edu}
\thanks{This work was supported in part by grants from the NSA Mathematical Sciences Program and the PSC-CUNY Research Award Program.}


\date{\today}

\begin{document}

\begin{abstract}
Let $f(x_1,x_2,\ldots,x_m) = u_1x_1+u_2 x_2+\cdots + u_mx_m$ be a linear form with positive integer coefficients, 
and let $N_f(k) = \min\{ |f(A)| : A \subseteq \Z \text{ and } |A|=k \}.$ 
A minimizing $k$-set for $f$ is a set $A$ such that $|A|=k$ and $|f(A)|  =  N_f(k).$   A finite sequence $(u_1, u_2,\ldots,u_m)$ of positive integers is called complete if $\left\{\sum_{j\in J}u_j : J \subseteq \{1,2, \ldots,m\} \right\} = \{0,1,2,\ldots, U\},$ where $U = \sum_{j=1}^mu_j.$   It is proved that if $f$ is an $m$-ary linear form whose coefficient sequence $(u_1,\ldots,u_m)$ is complete, then $N_f(k) = Uk-U+1$ and the minimizing $k$-sets are precisely the arithmetic progressions of length $k$.  Other extremal results on linear forms over finite sets of integers are obtained.
\end{abstract}

\maketitle

\section{Extremal functions for linear forms}

Let $m \geq 1$ and let $f:\Z^m \rightarrow \R$ be a real-valued function of $m$ integer variables.  For every finite set $A$ of integers, consider the set 
\[
f(A) = \{ f(a_1,a_2, \ldots , a_m) : a_1,a_2,\ldots,a_m \in A\}.
\]
Let $|A|$ denote the cardinality of the set $A$.  We define the functions
\[
N_f(k) =  \min\{ |f(A)| : A \subseteq \Z \text{ and } |A|=k \}
\]
and
\[
M_f(k) =  \max\{ |f(A)| : A \subseteq \Z \text{ and } |A|=k \}.
\]
A set $A$ with $|A|=k$ is called a $k$-set.  If $A$ is a $k$-set and $|f(A)|=N_f(k),$ then $A$ is called a \emph{minimizing $k$-set for $f$.} If $A$ is a $k$-set and $|f(A)|=M_f(k),$ then $A$ is called a \emph{maximizing $k$-set for $f$.}  An important  inverse problem in number theory is to compute the extremal functions $N_f(k)$ and $M_f(k)$, and to classify the minimizing and maximizing $k$-sets for $f$.

In this paper we study linear forms.  A classical example in additive number theory is the linear form  $f(x_1,x_2,\ldots,x_m)=x_1+x_2+\cdots + x_m.$   In this case, $N_f(k) = mk-m+1$ and the minimizing $k$-sets are the arithmetic progressions of length $k$ (Nathanson~\cite[Theorem 1.6]{nath96bb}).  Also, $M_f(k) = {k+m-1 \choose m} = k^{m}/m! + O\left(k^{m-1}\right)$ and the maximizing $k$-sets are sets $A$ of positive integers (called Sidon sets or $B_h$-sets) such that every integer has at most one representation as the sum of $h$ not necessarily distinct elements of $A$.
 
Let $\mathcal{LF}(m)$ denote the set of all $m$-ary linear forms
\[
f(x_1,\ldots,x_m) = u_1x_1+u_2 x_2+\cdots + u_mx_m
\]
with positive integer coefficients.   The function $f$ is called the \emph{linear form associated to the sequence $(u_1,\ldots,u_m).$}   Without loss of generality we can assume that 
\[
1 \leq u_1 \leq u_2 \leq \cdots \leq u_m 
\]
and
\[
\gcd(u_1,u_2,\ldots,u_m)=1.
\]
Let $\mathcal{LF}^*(m)$ denote the set of all $m$-ary linear forms with pairwise distinct coefficients, that is, with   
\[
1 \leq u_1 < u_2 < \cdots < u_m.
\]
We define the extremal functions 
\[
\mathcal{N}_m(k)  = \min\{  N_f(k) : f \in \mathcal{LF}(m) \}
\]
and
\[
\mathcal{N}^*_m(k)  = \min\{  N_f(k) : f \in \mathcal{LF}^*(m) \}.
\]
Since the only linear form with $m=1$ is $f(x_1)=u_1x_1,$ it follows that 
$\mathcal{N}_1(k)  = k$ for all $k \geq 1.$  For binary forms we have $\mathcal{N}_2(2) = 3$ and $\mathcal{N}^*_2(2)  = 4.$

Fix the integer $m \geq 1$.  For $k\geq 2$ and $f \in \mathcal{LF}(m)$, let $A$ be a $k$-set such that $N_f(k) = |f(A)|$ and let $a^{\prime} = \max(A)$ and $A' = A \setminus \{a^{\prime} \}.$   
Since $f(A) \supseteq f(A^{\prime})$ and $f(a^{\prime},\ldots,a^{\prime}) > \max(f(A^{\prime}))$, it follows that 
\[
N_f(k) = |f(A)| > |f(A^{\prime})| \geq N_f(k-1)
\]
and so $\mathcal{N}_m(k) > \mathcal{N}_m(k-1)$ and $\mathcal{N}^*_m(k) > \mathcal{N}^*_m(k-1)$ for all $k\geq 2.$  
Choosing a $(k-1)$-set $A'$ such that $M_f(k-1) = |f(A')|$ and an integer $a' > \max(A'),$ we define the $k$-set $A = A' \cup \{a'\}$.   
Since $f(A) \supseteq f(A^{\prime})$ and $f(a^{\prime},\ldots,a^{\prime}) > \max(f(A^{\prime}))$, we have 
\[
M_f(k-1) = |f(A^{\prime})| <  |f(A)| \leq M_f(k).
\]
Similarly, for fixed $k \geq 2,$ the extremal functions $\mathcal{N}_m(k)$ and $\mathcal{N}^*_m(k)$ are strictly increasing in $m$.

Denote the interval of integers $\{n\in \Z: x \leq n \leq y\}$ by $[x,y]$  .
Given a finite sequence of integers $\mathcal{U} = (u_1, u_2, \ldots, u_m)$, we define the set of subset sums 
\[
S(\mathcal{U}) 
= \left\{ \sum_{j\in J} u_j : J \subseteq [1,m] \right\}.
\]
Let $U = \sum_{j=1}^m u_j.$   Then 
\[
\{ 0,U \} \subseteq S(\mathcal{U}) \subseteq [0,U]
\]
and $n \in S(\mathcal{U})$ if and only if $U-n \in S(\mathcal{U})$.  
The sequence $\mathcal{U}$ is called \emph{complete} if  $S(\mathcal{U}) = [0,U]$.  For example, the sequence $(1,2,3,\ldots,m)$ is complete for all $m\geq 1.$   The sequence $(1,1,3)$ is complete, but the sequence $(1,3)$ is not complete.    
(This is the finite analogue of an infinite complete sequence, which is a sequence $\mathcal{U}$ of positive integers such that $S(\mathcal{U})$ 
contains all sufficiently large integers (cf. Szemer\' edi-Vu~\cite{szem-vu06}).)

The sequence $\mathcal{U}$ has \emph{distinct subset sums} if $|S(\mathcal{U})| = 2^m,$ that is,  if the conditions $I,J \subseteq \{1,2,\ldots,m\}$ and $\sum_{i\in I} u_i = \sum_{j\in J} u_j$ imply that $I=J$.  For example, the sequence $(1,g,g^2,\ldots,g^{k-1})$ has distinct subset sums for every $g\geq 2.$

If $\mathcal{U} = (u_1,  \ldots, u_m)$ is an increasing sequence of positive integers and $f(x_1,\ldots,x_m) = \sum_{j=1}^m u_j x_j,$ then 
\[
f(1,\ldots,1) = \sum_{j=1}^m u_j = U
\]
and
\[
f\left( \{ 0,1 \} \right)  = S(\mathcal{U}).
\]
In particular, $\mathcal{U}$ is complete if and only if $f\left( \{ 0,1 \} \right) = [0,U],$ and $\mathcal{U}$ has discrete subset sums if and only if $\left| f\left( \{ 0,1 \} \right) \right| = 2^m=M_f(2).$

For a finite set $A$ of integers  and for integers $c \neq 0$ and $d$, we define the affine transformation 
\[
c\ast A+d = \{ca+d : a\in A \}.
\]
If $f(x_1,\ldots,x_m)=u_1x_1+\cdots + u_mx_m\in \mathcal{LF}(m)$ with $U = u_1+\cdots + u_m,$ then 
\[
f(c\ast A + d) = c\ast f(A) + dU
\]
and 
\[
|f(c\ast A + d)| = |f(A)|
\]
for integers $c \neq 0$ and $d$.  Thus, the function $|f(A)|$ is an affine invariant of $A$ (cf. Nathanson~\cite{nath07a}).

The study of inverse problems for $m$-ary forms is related to the paper~\cite{nath07d}, which initiated the comparative study of binary linear forms.

\section{A lower bound for $m$-ary linear forms}
The following result is elementary but fundamental.

\bl   \label{LFB:lemma:general-main}
Let $f:\Z^m \rightarrow \R$ be a real-valued function of $m$ integer variables.  Let  $g:\Z \rightarrow \R$ be a strictly increasing function such that 
\beq  \label{LFB:condition-i}
f([a,b]) \subseteq [g(a),g(b)] \qquad\text{for all integers $a < b$.}
\eeq
Let $\ell$ and $\lambda$ be positive integers with $\ell \geq 2$ such that
\beq  \label{LFB:condition-ii}
N_f(\ell) \geq \lambda.
\eeq
Let $A = \{a_i\}_{i=0}^{k-1}$ be a set of $k$ integers with $a_{i-1}<a_i$ for $i=1,\ldots,k-1.$  Let $k-1=q(\ell -1)+r,$ where $0 \leq r \leq \ell -2.$  Define 
\[
\mu(A) = \left| f\left( \{ a_{k-r-1},a_{k-r},a_{k-r+1},\ldots,a_{k-1} \}\right) \right|.
\]
Then
\beq  \label{LFB:main-ineq-i}
|f(A)|  \geq (\lambda -1)\left(\frac{k-r-1}{\ell -1} \right) +  \mu(A)
\eeq
and
\beq   \label{LFB:main-ineq-ii}
N_f(k) \geq \left(\frac{\lambda -1}{\ell -1}\right) k - \lambda +2
\eeq
for every positive integer $k$.
\el

\begin{proof}
Dividing $k$ by $\ell -1,$ we have
$k - 1 = q(\ell-1)+r,$
where $0 \leq r \leq \ell-2.$  Then $k \leq (q+1)(\ell -1).$
Let $A = \{a_0,a_1, \ldots,a_{k-1} \}$ be a set of $k$ integers, where 
\[
a_0 < a_1 < \cdots < a_{k-1}.
\]
For $j=0,1,\ldots,q-1,$ we define the sets
\[
A_j =   \left\{ a_i : i \in [j(\ell -1), (j+1)(\ell -1)] \right\}.
\] 
Let $A_q$ be the set 
\[
A_q =  \left\{ a_i : i \in [ q(\ell -1) ,  q(\ell -1) + r] \right\} 
=   \left\{a_{k-r-1},a_{k-r},a_{k-r+1},\ldots,a_{k-1} \right\}
\]
Then
\[
A_{j-1} \cap A_j = \{j(\ell -1)\}
\]
for $j=1,\ldots,q,$ and 
\[
A = \bigcup_{j=0}^{q} A_j.
\]
Since 
\[
\max\left( A_{j-1} \right) = a_{j(\ell -1)} = \min\left( A_{j} \right)
\]
for $j=1,\ldots,q,$ and since the function $g$ is strictly increasing, condition~~\eqref{LFB:condition-i} implies that $f(A_{j'}) \cap f(A_{j}) = \emptyset $ if  $j - j' >1$, and 
\[
f\left(A_{j-1} \right) \cap f\left( A_{j} \right) =  \{g( j(\ell -1)) \}  = \{ f(j(\ell -1),\ldots,j(\ell -1)) \}
\]
for $j=1,\ldots,q$.  
Note that $\mu(A) = |f(A_q)| = 1$ if $r=0.$

By condition~\eqref{LFB:condition-ii}, we have $|f(A_j)| \geq \lambda$ for $j=0,1,\ldots,q-1,$ and so 
\begin{align*}
|f(A)| & = \left| f\left( \bigcup_{j=0}^{q} A_j \right) \right|  \\
& \geq \sum_{j=0}^{q} |f(A_j)| - q \\
& \geq \lambda  q- q  + \mu(A) \\
& = (\lambda -1)\left(\frac{k-r-1}{\ell -1} \right) + \mu(A)\\
& \geq \left(\frac{\lambda -1}{\ell -1}\right) k - \lambda + 2.
\end{align*}
The observation that the last inequality is independent of the set $A$ completes the proof.
\end{proof}

\bl   \label{LFB:lemma:main}
Let $f \in \mathcal{LF}(m).$    If $\ell$ and $\lambda$ are positive integers with $\ell \geq 2$ such that $N_f(\ell) \geq \lambda,$ then
\[
N_f(k) \geq \left(\frac{\lambda -1}{\ell -1}\right) k - \lambda +2
\]
for every positive integer $k$.
\el

\begin{proof}
Let $f(x_1,\ldots,x_m) = \sum_{j=1}^m u_j x_j.$  
For $U = \sum_{j=1}^m u_j,$  we define the strictly increasing function $g(x) = Ux.$  If  $a \leq x_j \leq b$ for $j=1,\ldots,m,$ then
\[
g(a) = Ua \leq f(x_1,\ldots,x_m) \leq Ub = g(b)
\]
and so
\[
f([a,b]) \subseteq [g(a),g(b)]
\]
for all integers $a < b$. 
The result follows from Lemma~\ref{LFB:lemma:general-main}.
\end{proof}

\bt   \label{LFB:theorem:m-ary}
For all positive integers $m$ and $k,$
\[
\mathcal{N}^*_m(k) = \left(  \frac{m^2+m}{2} \right) k - \left( \frac{m^2+m-2}{2} \right).
\]
\et

\begin{proof}
Let $f \in \mathcal{LF}^*(m)$.  Then $f(x_1,x_2,\ldots,x_m)=u_1x_1+u_2x_2+ \cdots + u_mx_m$ with $1 \leq u_1 < u_2 < \cdots < u_m.$
For integers $a<b$ and $i=0,1,\ldots,m$, we define the integer 
\begin{align*}
s_i & = f\left(\underbrace{a,\ldots,a}_{m-i \text{ terms}},\underbrace{b,\ldots,b}_{i \text{ terms}} \right) \\
& = \left(u_1+\cdots + u_{m-i}\right)a + \left(u_{m-i+1}+\cdots + u_m\right)b \\
& \in f(A).
\end{align*}
Then 
\[
s_0 < s_1 < \cdots < s_m.
\]
For $i = 0,1,\ldots, m-2$ and $j=0,1,\ldots,m-i,$ the integer
\begin{align*}
t_{i,j} & =  f\left( \underbrace{a,\ldots,a}_{j-1 \text{ terms}},b,\underbrace{a,\ldots,a}_{m-i-j\text{ terms}},\underbrace{b,\ldots,b}_{i \text{ terms}} \right) \\
& = \left(u_1+\cdots + u_{j-1}\right)a + u_j b + \left(u_{j+1}+\cdots + u_{m-i}\right)a +  \left(u_{m-i+1}+\cdots + u_{m}\right)b \\
& \in f(A)
\end{align*}
satisfies
\[
s_i  = t_{i,0} < t_{i,1} < \cdots < t_{i,m-i-1} < t_{i,m-i} = s_{i+1},
\]
It follows that 
\[
N_f(2) \geq (m+1) + \sum_{i=0}^{m-2} (m-i-1) = \frac{m^2+m+2}{2}.
\]
Applying Lemma~\ref{LFB:lemma:main} with $\ell = 2$ and $\lambda = (m^2+m+2)/2,$ we obtain
\[
\mathcal{N}^*_m(k) \geq \left(  \frac{m^2+m}{2} \right) k - \left( \frac{m^2+m-2}{2} \right).
\]

To prove that this lower bound is best possible, we consider the linear form 
\[
f(x_1,\ldots,x_m) = x_1 + 2x_2 + \cdots + ix_i + \cdots + mx_m \in \mathcal{LF}^*(m)
\]
and the finite set 
\[
A = \{ 0,1,\ldots, k-1\}.
\]
Then 
\[
f(A) = \left[ 0,  \frac{m(m+1)(k-1)}{2} \right]
\]
and so
\[
|f(A)| = \left(  \frac{m^2+m}{2} \right) k - \left( \frac{m^2+m-2}{2} \right) = \mathcal{N}^*_m(k).
\]
This completes the proof.
\end{proof}

\section{A lower bound for binary and ternary linear forms}

\bt     \label{LFB:theorem:binary}
Let $f(x_1,x_2)=u_1x_1 +u_2x_2 \in \mathcal{LF}(2),$ where $1 \leq u_1 <  u_2 $ and  $\gcd(u_1,u_2)=1$.
\benum
\item[(i)]
If $f(x_1,x_2)= x_1 + x_2,$ then $N_f(k)=2k-1.$
\item[(ii)]
If $f(x_1,x_2)= x_1 + 2x_2,$ then $N_f(k)= 3k-2.$
\item[(iii)]
If $f(x_1,x_2)\neq x_1 + x_2$ or $x_1 + 2x_2,$ then 
\[
N_f(k) \geq \left[\frac{7k-5}{2}\right].
\]
\eenum
\et

\begin{proof}
Let $|A|=k.$  If $f(x_1,x_2)= x_1 + x_2,$ then $|f(A)| \geq 2k-1$ and $f([0,k-1]) = [0,2k-2]$, hence $|f([0,k-1])| = 2k-1.$   

If $f(x_1,x_2)= x_1 + 2x_2,$ then $|f(A)| \geq 3k-2$ by Theorem~\ref{LFB:theorem:m-ary}.  Moreover,  $f([0,k-1]) = [0,3k-3]$ and so $|f([0,k-1] = 3k-2.$   

If $f(x_1,x_2)=u_1x_1 +u_2x_2 \in \mathcal{LF}(2)$ and $f(x_1,x_2)\neq x_1 + x_2$ or $x_1 + 2x_2,$ then $u_2 \geq 3.$  We shall prove that  $N_f(3) = 8$ or $9$.
We use the fact that the quadratic form $u_1^2 + u_1u_2-u_2^2\neq 0$ for all nonzero integers $u_1$ and $u_2.$  

Let $A=\{a,b,c\},$ where $a < b < c.$    Then $|f(A)| \leq 9.$  We have the following strictly increasing sequence of seven elements of $f(A)$:
\begin{align*}
u_1a + u_2 a & <  u_1b + u_2a  < u_1a + u_2b  < u_1b + u_2b  \\
& < u_1c + u_2  b < u_1b + u_2c  < u_1c + u_2 c
\end{align*}
and so $ |f(A)| \geq 7.$   There is another strictly increasing sequence of four elements of $f(A)$:
\[
u_1b + u_2a  < u_1c + u_2a  < u_1a + u_2 c < u_1b + u_2c.
\]
If $|f(A)|=7,$ then 
\beq  \label{LFB:setincl}
\{ u_1c + u_2a, u_1a + u_2 c\} \subseteq \{ u_1a + u_2b  , u_1b + u_2b,  u_1c + u_2  b \}.
\eeq
This is possible in only three ways.  In the first case, we have 
\begin{align*}
 u_1c + u_2a  & =  u_1a + u_2b   \\
 u_1a + u_2 c & =  u_1b + u_2b.
\end{align*}
Eliminating $a$ from these equations, we obtain $(u_1^2 + u_1u_2-u_2^2)(c-b) =0,$ which is false.

In the second case,
\begin{align*}
 u_1a + u_2c  & =  u_1b + u_2a   \\
 u_1c + u_2 a & =  u_1b + u_2  c.
\end{align*}
Eliminating $b$ from these equations, we obtain $(u_2-2u_1)(c-a)=0$ and so $2u_1=u_2.$  Since $\gcd(u_1,u_2)=1,$ it follows that $u_1=1$ and $u_2=2,$ which is also false.

In the third case,
\begin{align*}
 u_1c + u_2a  & =  u_1b + u_2b   \\
 u_1a + u_2c & =  u_1c + u_2  b.
\end{align*}
Eliminating $a$ from these equations, we again obtain $(u_1^2 + u_1u_2-u_2^2)(c - b) = 0,$ which is false.  
It follows that~\eqref{LFB:setincl} is impossible, and so $|f(A)| \geq 8.$  
  Applying Lemma~\ref{LFB:lemma:main} with $\ell=3$ and $\lambda = 8,$ we obtain 
\[
N_f(k) \geq \frac{7k}{2} -6.
\]

We can improve the constant term by using the more precise inequality~\eqref{LFB:main-ineq-i} in Lemma~\ref{LFB:lemma:general-main}.  If $r=0$, then $\mu(A)=1$ and  
\[
|f(A)|  \geq 7\left(\frac{k-1}{2} \right) + 1 
= \frac{7k-5}{2}.
\]
If $r=1$, then $\mu(A) = N_f(2) = 4$ and  
\[
|f(A)|  \geq 7\left(\frac{k-2}{2} \right) + 4 
= \frac{7k-6}{2}.
\]
This completes the proof.
\end{proof}

\bl   \label{LFB:lemma:ternary}
Let $f(x_1,x_2,x_3)=u_1x_1+u_2x_2+u_3x_3 \in \mathcal{LF}(3)$ with 
$1 \leq u_1 \leq u_2 \leq u_3$ and $\gcd(u_1,u_2,u_3)=1.$ 
If $f \in \mathcal{LF}^*(3)$, then $N_f(2)=7$ or $8,$ and $N_f(2)=8$ if and only if $u_1 + u_2 \neq u_3.$   Also, $\mathcal{N}^*_3(2) = 7.$

Let $f \in \mathcal{LF}(3) \setminus \mathcal{LF}^*(3).$
\benum
\item[(i)]
If $u_1=u_2=u_3=1$, then $N_f(2)=4$.
\item[(ii)]
If $u_1=u_2$ and  $u_3=2u_1$, then $N_f(2)=5$.
\item[(iii)]
If $u_1=u_2$ and  $u_3 \neq 2u_1$, then $N_f(2)=6$.
\item[(iv)]
If $u_1<  u_2=u_3$, then $N_f(2)=6.$
\eenum
\el

\begin{proof}
Let $f(x_1,x_2,x_3)=u_1x_1+u_2x_2+u_3x_3$, where  $1 \leq u_1  < u_2 < u_3$  Then
\begin{align*}
u_1a+u_2a +u_3 a & < u_1b+u_2 a + u_3 a < u_1a + u_2b+u_3a \\
& < u_1a+u_2a + u_3b < u_1b+u_2a+u_3b \\
& < u_1a+u_2b+u_3b < u_1b+u_2b+u_3b.
\end{align*}
These inequalities account for seven of the at most eight elements of the set $f(A).$  The remaining element is $f(a,b,b) = u_1b+u_2b+u_3a.$ 
Since
\[
u_1a+u_2b +u_3a < u_1b+u_2b+u_3a < u_1b+u_2a+u_3b
\]
it follows that $N_f(2)=7$ if and only if
$u_1a +u_2a +u_3b = u_1b +u_2b+u_3a.$
This is equivalent to $(u_1+u_2-u_3)(b-a)=0$ or $u_1 + u_2 = u_3.$
It follows that $\mathcal{N}^*_3(2) = N_f(2) = 7$ if and only if $u_1 + u_2 = u_3.$

Identities~(i)-(iv) are straightforward calculations.
\end{proof}

\bt
Let $f(x_1,x_2,x_3)=u_1x_1+u_2x_2+u_3x_3 \in \mathcal{LF}^*(3)$ with 
$1 < u_1 < u_2 < u_3$ and $\gcd(u_1,u_2,u_3)=1.$  Then $N_k(f) \geq 6k-5.$  If $f \in \mathcal{LF}^*(3)$ and $u_1 + u_2 \neq u_3,$  then $N_f(k) \geq 7k-6.$
\et

\begin{proof}
Applying Theorem~\ref{LFB:theorem:m-ary} with $m=3$ gives $N_k(f) \geq 6k-5.$   By Lemma~\ref{LFB:lemma:ternary}, if $f \in \mathcal{LF}^*(3)$ and $u_1 + u_2 \neq u_3,$  then $N_f(2)=8.$   Applying Lemma~\ref{LFB:lemma:main} with $\ell = 2$ and $\lambda = 8$ gives $N_f(k) \geq 7k-6.$
\end{proof}

Note that an increasing sequence $(u_1, u_2, u_3)$ has distinct subset sums if and only if it is strictly increasing and $u_1+u_2 \neq u_3.$

\section{An inverse problem for linear forms}
Let $f$  be a linear form in $m$ variables with positive integral coefficients.  The inverse problem for $f$ is to determine the $k$-minimizing sets for $f$, that is, to describe the structure of a $k$-set $A$ such that $|f(A)| = N_f(k).$  For example, if $f(x_1,\ldots,x_m)=x_1+ \cdots + x_m,$  then $N_f(k)= mk-m+1,$ and $N_f(A) = mk-m+1$ if and only if $A$ is an arithmetic progression of length $k$  (Nathanson~\cite[Theorem 1.6]{nath96bb}).  If $f(x_1,x_2)=x_1+ 2x_2,$  then Cilleruelo,  Silva, and Vinuesa~\cite{cill-silv-vinu07}  proved that $N_f(k)=3k-2,$ and $N_f(A)=3k-2$ if and only if $A$ is an arithmetic progression.  This result generalizes to all $m$-ary forms whose coefficient sequence is complete.

\bt
Let $\mathcal{U} = (u_1,\ldots,u_m)$ be a complete increasing sequence of positive integers with $U = \sum_{j=1}^m u_j.$   Consider the linear form  $f(x_1,\ldots,x_m) = u_1x_1 + \cdots + u_m x_m$.  Then $N_f(k) = Uk-U+1$, and the set $A$ is a minimizing $k$-set for $f$ if and only if $A$ is an arithmetic progression of length $k$.
\et

\begin{proof}
Since $\mathcal{U}$ is complete, it follows that
for any integers $a$ and $b$ with $a < b$ we have 
\begin{align*}
f(\{a,b \}) 
& = \left\{ \left( \sum_{i\in  [1,m]\setminus I} u_i \right)a 
            +  \left( \sum_{i\in I}u_i \right)b  : I \subseteq [1,m] \right\} \\
& = \left\{ \left( U-\ell  \right)a  + \ell b : \ell = 0,1,\ldots,U \right\} \\
& = \left\{ Ua + \ell (b-a) : \ell = 0,1,\ldots,U \right\}.
\end{align*}
Since $f(\{ i-1, i\}) = [U(i-1),Ui]$, it follows that
\begin{align*}
[0,U(k-1)] & = \bigcup_{i=1}^{k-1} [U(i-1),Ui] = \bigcup_{i=1}^{k-1} f(\{i-1,i\}) \\
& \subseteq f([0,k-1]) \subseteq [0,U(k-1)].
\end{align*}
Then $f([0,k-1]) = [0,U(k-1)]$ and $N_f(k) \leq |f([0,k-1])|=Uk-U+1.$  

Applying Lemma~\ref{LFB:lemma:main} with $\ell = 2$ and $\lambda =  U+1,$ we obtain the lower bound $|f(A)| \geq Uk-U+1,$ 
and so $N_f(k)=Uk-U+1.$
Since $|f([0,k-1])|=Uk-U+1$ and $|f(A)|$ is an affine invariant of $A$, it follows that $|f(A)| = Uk-k+1$ for every arithmetic progression $A$ of length $k$.

Conversely, let $A = \{ a_0, a_1,\ldots,a_{k-1} \}$ be a minimizing $k$-set for $f$ with $a_0 < a_1 < \cdots < a_{k-1}.$   Since $(u_1,\ldots,u_m)$ is a complete sequence, 
\[
f(\{a_{i-1}, a_i \})  = \left\{ \left( U-\ell  \right)a_{i-1}  + \ell a_i : \ell = 0,1,\ldots,U \right\}.
\]
For $i=1,\ldots, k-2$ we have the inequalities 
\begin{align*}
 \boxed{(U-1)a_{i-1}+a_i } & <  (U-2)a_{i-1}+ 2a_i < \cdots  <  a_{i-1}+ (U-1)a_i < Ua_i   \\
  & < (U-1)a_i+a_{i+1} <  \cdots 
 < 2a_i+ (U-2)a_{i+1}  \\
 & < \boxed{a_i+ (U-1)a_{i+1}} < Ua_{i+1}.
\end{align*}
Since $|f(A)| = Uk-U+1$, it follows that 
\beq    \label{LFB:structure}
f(A) = \bigcup_{i=1}^{k-1} f(\{ a_{i-1}, a_i\}) 
= \bigcup_{i=1}^{k-1} \left\{ (U-k)a_{i-1}  + ka_i : k=0,1,\ldots,U \right\}.
\eeq
We also have
\begin{align*}
\boxed{ (U-1)a_{i-1} + a_i  }& <   (U-1)a_{i-1}+a_{i+1} < (U-2)a_{i-1}+ 2a_{i+1} < \cdots \\
& < 2a_{i-1}+ (U-2)a_{i+1} <  a_{i-1}+ (U-1)a_{i+1}   \\
& < \boxed{ a_i+ (U-1)a_{i+1} }.
\end{align*}
Equation~\eqref{LFB:structure} implies that 
\begin{align}  \label{LFB:inclusion}
\{ (U-k) & a_{i-1}  +ka_{i+1}: k=1,\ldots,U-1 \}   \nonumber  \\
\subseteq &
\{ (U-k)a_{i-1}+ka_i: k=2,\ldots,U \}   \\
& \cup \{ (U-k)a_i + ka_{i+1}: k=1,\ldots,U-2 \}.  \nonumber
\end{align}

We want to prove that $A$ is an arithmetic progression.  If not, then $a_{i-1}+a_{i+1} \neq 2a_i$ for some $i \in [1,k-2].$  
It follows that for all $k \in [1,U/2]$ we have
\beq   \label{LFB:inverse1}
(U-k)a_{i-1}+ka_{i+1} \neq (U-2k)a_{i-1}+2ka_i
\eeq
and
\beq   \label{LFB:inverse2}
ka_{i-1} + (U-k)a_{i+1} \neq 2ka_i+ (U - 2k) a_{i+1}.
\eeq

Let $U' = U/2$ if $U$ is even and $U' = (U-1)/2$ if $U$ is odd.  
Set inclusion~\eqref{LFB:inclusion} implies that 
\[
(U-1)a_{i-1}+a_{i+1} \geq (U-2)a_{i-1}+2a_i.
\]
Suppose that 
\[
(U-k)a_{i-1} + ka_{i+1} \geq (U-2k)a_{i-1} + 2k a_i
\]
for some $k \in [1, U'-1].$  
We deduce from inequality~\eqref{LFB:inverse1} that 
\[
(U-k)a_{i-1} + ka_{i+1} > (U-2k)a_{i-1} + 2k a_i
\]
and so, again by~\eqref{LFB:inverse1},
\[
(U-k)a_{i-1}+ ka_{i+1} \geq (U-2k-1)a_{i-1}+(2k+1)a_i.
\]
It follows again from~\eqref{LFB:inclusion} that 
\[
(U-(k+1))a_{i-1}+ (k+1)a_{i+1} \geq (U-2(k+1))a_{i-1}+ 2(k+1)a_i.
\]
Continuing inductively, we obtain 
\beq   \label{LFB:Kineq}
(U-U')a_{i-1}+U'a_{i+1} \geq (U-2U')a_{i-1}+2U'a_i
\eeq
and so
\[
(U-U')a_{i-1}+U'a_{i+1} > (U-2U')a_{i-1}+2U'a_i.
\]
If $U = 2U'$ is even, this inequality can be rewritten as 
\[
U'a_{i-1}+U'a_{i+1} \geq Ua_i.
\]
If $U = 2U'+1$ is odd, inequality~\eqref{LFB:Kineq} becomes
\[
(U'+1)a_{i-1}+U'a_{i+1} \geq a_{i-1} + (U-1) a_i.
\]
Inequality~\eqref{LFB:inverse1} and set inclusion~\eqref{LFB:inclusion} imply that 
\[
(U'+1)a_{i-1}+U'a_{i+1} \geq U a_i.
\]
Therefore,
\[
U'a_{i-1}+ (U' + 1) a_{i+1} \geq (U-1) a_i + a_{i+1}.
\]
In both cases we have
\beq   \label{LFB:downineq}
ka_{i-1}+ (U-k) a_{i+1} \geq 2k a_i + (U - 2k) a_{i+1}
\eeq
for $k = U'.$

Suppose that~\eqref{LFB:downineq} holds for some $k\in [2,U'].$  
Inequality~\eqref{LFB:inverse2} and set inclusion~\eqref{LFB:inclusion} imply that 
\[
ka_{i-1}+ (U-k) a_{i+1} \geq (2k-1) a_i + (U - (2k-1)) a_{i+1}.
\]
Therefore,
\[
(k-1)a_{i-1}+ (U-(k-1)) a_{i+1} \geq 2(k-1) a_i + (U - 2(k-1)) a_{i+1}.
\]
Continuing downward inductively, we obtain
\[
a_{i-1}+ (U-1) a_{i+1} \geq 2 a_i + (U - 2) a_{i+1}.
\]
Since $a_{i-1}+ (U-1) a_{i+1} < a_i+ (U-1) a_{i+1},$ it follows that 
$a_{i-1}+ (U-1) a_{i+1} =  2 a_i + (U - 2) a_{i+1},$ which implies that $a_{i-1}+ a_{i+1} =  2 a_i.$  This is a contradiction.  Therefore, the minimizing $k$-set $A$ is an arithmetic progression.  This completes the proof.
\end{proof}

\section{An upper bound for linear forms}
We record here some simple estimates for the maximal function $M_f(k).$ 

\bt
For all $m$-ary linear forms $f \in \mathcal{LF}(m)$ and all positive integers $k,$
\beq  \label{LFB:1}
k^m \geq M_f(k) \geq {k \choose m}.
\eeq
If $f \in \mathcal{LF}^*(m),$ then 
\beq  \label{LFB:2}
M_f(k) \geq k(k-1)\cdots (k-m+1).
\eeq
If $U = \{u_1, u_2, \ldots,u_m\}$ is an increasing sequence of positive integers with distinct subset sums, and $f(x_1,\ldots,x_m) = u_1 x_1+u_2 x_2+\cdots + u_m x_m \in \mathcal{LF}(m)$, then
\beq  \label{LFB:3}
M_f(k) = k^m.
\eeq
\et

\begin{proof}
Let $f(x_1,\ldots,x_m) = u_1 x_1+u_2 x_2+\cdots + u_m x_m \in \mathcal{LF}(m)$.  The upper bound for $M_f(k)$ comes from counting the number of $m$-tuples of a $k$-element set.    
To obtain the lower bound in~\eqref{LFB:1}, choose $g > mu_m$ and  let $A=\{1,g ,g^2,\ldots,g^{k-1} \}.$   If $(r_1,\ldots,r_m)$ and $(s_1,\ldots,s_m)$ are $m$-tuples of elements of $[0,k-1]$ such that 
$|\{r_1,\ldots,r_m\}| = |\{s_1,\ldots,s_m\}| = m$ and the $k$-sets 
$\{r_1,\ldots,r_m\} \neq \{s_1,\ldots,s_m\}$ are distinct, then the uniqueness of the $g$-adic representations of the positive integers implies that $f(g^{r_1},\ldots,g^{r_m}) \neq f(g^{s_1},\ldots, g^{s_m})$.  This proves~\eqref{LFB:1}. 

If $f\in \mathcal{LF}^*(m),$  then the coefficients $u_1,\ldots,u_m$ are pairwise distinct. If $(r_1,\ldots,r_m)$ and $(s_1,\ldots,s_m)$ are $m$-tuples of elements of $[0,k-1]$ such that 
$|\{r_1,\ldots,r_m\}| = |\{s_1,\ldots,s_m\}| = m$ and the $m$-tuples
$(r_1,\ldots,r_m)$ and $(s_1,\ldots,s_m)$ are distinct, then the uniqueness of the $g$-adic representations of the positive integers implies that $f(g^{r_1},\ldots,g^{r_m}) \neq f(g^{s_1},\ldots, g^{s_m})$.  This proves~\eqref{LFB:2}. 

Finally, suppose that the sequence $(u_1, u_2, \ldots,u_m)$ has distinct subset sums.  Let $(r_1,\ldots, r_m)$ and  $(s_1,\ldots, s_m)$ be $m$-tuples of elements of $[0,k-1]$ such that  $f(g^{r_1},\ldots, g^{r_m}) = f(g^{s_1},\ldots, g^{s_m}).$   For $d \in [0,k-1],$ let $I_d = \{i \in [1,m]: r_i =d\}$ and $J_d = \{j \in [1,m] : s_j=d\}$.  Then
\[
f(g^{r_1},\ldots, g^{r_m}) = \sum_{d=0}^{k-1} \left(\sum_{i \in I_d} u_i\right)g^d
\]
and
\[
f(g^{s_1},\ldots, g^{s_m}) = \sum_{d=0}^{k-1} \left(\sum_{j\in J_d} u_i\right)g^d.
\]
Since
\[
\max\left(  \sum_{i\in I_d} u_i , \sum_{j\in J_d} u_i \right) \leq  mu_m < g
\]
it follows that if $f(g^{r_1},\ldots, g^{r_m}) = f(g^{s_1},\ldots, g^{s_m}),$ then $I_d = J_d$ for all $d$ and so $r_i=s_i$ for $i=1,\ldots,m.$    Therefore, $|f(A)|=k^m.$
\end{proof}

\section{Open problems}
\benum
\item
The minimizing $k$-sets for linear forms associated to complete sequences are precisely the arithmetic progressions of length $k$.   Classify all linear forms $f(x_1,\ldots, x_m)$ with the property that the only minimzing $k$-sets are arithmetic progressions.  In particular, if $f(x_1,\ldots, x_m) = u_1x_1+\cdots + u_m x_m$ is a linear form whose minimizing $k$-sets are arithmetic progressions, then is the sequence $(u_1,\ldots,u_m)$ complete?

\item
Let $f(x_1,\ldots, x_m) = u_1x_1+\cdots + u_m x_m$ be a linear form with $U = \sum_{j=1}^m u_j.$  Is the sequence $(u_1,\ldots,u_m)$ complete if $N_f(k) = Uk-U+1$?

\item
There is no reason to consider only linear forms. Let $f(x_1,\ldots,x_m)$ be a polynomial with integer coefficients.  The set $A$ is a \emph{minimizing $k$-set for $f$} if $|f(A)| = N_f(k).$  Compute $N_f(k)$ and determine the minimizing $k$-sets for $f$.

\item
Let $s(x_1,x_2)=x_1+x_2.$  The Freiman philosophy of inverse problems in additive number theory is to deduce structural  information about a finite set $A$ of integers if the sumset $s(A) = A+A$ is small (cf. Freiman~\cite{frei73}).  Analogously, a natural inverse problem for linear forms and, more generally, arbitrary integer-valued polynomials in $m$ variables, is to deduce information about the finite sets $A$ of integers such that $|f(A)|-N_f(A)$ is small.

\item
For $f\in \mathcal{LF}(m)$, define the set  
\[
\mathcal{E}_f(k) = \{|f(A)| : A \subseteq \Z \text{ and } |A| = k \}.
\]  
By definition, $\min\left( \mathcal{E}_f(k) \right) = N_f(k)$ and 
$\max\left( \mathcal{E}_f(k)\right) = M_f(k).$  
For example, if $f\in \mathcal{LF}(2)$, then $\mathcal{E}_f(2) = [3,4]$, and, by Lemma~\ref{LFB:lemma:ternary}, $\mathcal{E}_f(3) = [4,8]$.
When is the set $\mathcal{E}_f(k)$ an interval of integers?  
For every linear form $f$ and $e \in  \mathcal{E}_f(k)$, let $\mathcal{A}_f(e)$ be the set of all $k$-sets $A$ of integers such that $|f(A)| = e.$  Then  $\{ \mathcal{A}_f(e) \}_{e\in \mathcal{E}_f(k)}$ is a partition of the $k$-sets of integers.  Can one classify the sets in this partition?
 There are many such questions.

 \eenum

\emph{Acknowledgement.}
I wish to thank Manuel Silva for several useful conversations about   binary linear forms.

\def\cprime{$'$} \def\cprime{$'$} \def\cprime{$'$}
\providecommand{\bysame}{\leavevmode\hbox to3em{\hrulefill}\thinspace}
\providecommand{\MR}{\relax\ifhmode\unskip\space\fi MR }
\providecommand{\MRhref}[2]{%
  \href{http://www.ams.org/mathscinet-getitem?mr=#1}{#2}
}
\providecommand{\href}[2]{#2}

\end{document}